\documentclass{article}
\usepackage{amsfonts,amssymb,amsmath,amsthm}

\newtheorem{theorem}{Theorem}

\newcommand{\Prob}{\mathbb{P}}
\newcommand{\expct}{\mathbb{E}}
\newcommand{\diffd}{\mathrm{d}}
\newcommand{\R}{\mathbb R}

\newcommand{\eps}{\varepsilon}

\newcommand{\F}{{\mathcal{F}}}

\newcommand{\E}{\mathbb{E}}

\begin{document}

\title{\bf The almost-sure population growth rate in branching Brownian motion with a quadratic breeding potential}
\author{J. Berestycki\thanks{%
Laboratoire de Probabilit\'{e}s et Mod\`{e}les Al\'{e}atoires, Universit\'{e} Pierre et Marie Curie (UPMC Paris VI),  CNRS UMR 7599, 175 rue du Chevaleret, 75013 Paris, France. email:\texttt{julien.berestycki@upmc.fr.}
}\ , 
\'{E}. Brunet\thanks{%
Laboratoire de Physique Statistique, \'{E}cole Normale Sup\'{e}rieure, UPMC
Universit\'{e} Paris 6, Universit\'{e} Paris Diderot, CNRS, 24 rue Lhomond, 75005
Paris, France. email: \texttt{eric.brunet@lps.ens.fr.}
}\ ,
J. W. Harris\thanks{%
Department of Mathematics, University of Bristol, University Walk,
Bristol, BS8 1TW, U.K. email: \texttt{john.harris@bristol.ac.uk.} Supported by the Heilbronn Institute for Mathematical Research.} \ and S. C. Harris%
\thanks{%
Department of Mathematical Sciences, University of Bath, Claverton Down,
Bath, BA2 7AY, U.K. email: \texttt{s.c.harris@bath.ac.uk}} }
\maketitle

\begin{abstract}
In this note we consider a branching Brownian motion (BBM) on $\R$ in which a particle at spatial position $y$ splits into two at rate $\beta y^2$, where $\beta>0$ is a constant. This is a critical breeding rate for BBM in the sense that the expected population size blows up in finite time while the population size remains finite, almost surely, for all time. We find an asymptotic for the almost sure rate of growth of the population.
\medskip


\end{abstract}

\section{Introduction}\label{sec:intro}
We consider a branching Brownian motion with a quadratic
breeding potential. Each particle diffuses as a driftless
Brownian motion on $\R$, and splits into two particles at rate $\beta y^2$,
where $\beta>0$ and $y$ is the spatial position of the particle. We let $N_t$ be the set of particles alive at time $t$, and then, for
each $u\in N_t$, $Y_u(t)$ is the spatial position of particle $u$
at time $t$ (and for $0\le s < t, Y_u(s)$ is the spatial position of the unique ancestor of $u$ alive at time $s$). We will call this process the $(\beta
y^2;\R)$-BBM. 

It is known that
quadratic breeding is a critical rate for population explosions.
If the breeding rate were instead $\beta |y|^p$ for $p>2$, the
population size would almost surely explode in a finite time. However,
for the $(\beta y^2;\R)$-BBM the \emph{expected} number of
particles blows up in a finite time while the total number of
particles alive remains finite almost surely, for all time. For $p\in[0,2)$ the expected population size remains finite for all time. See It\^{o} and
McKean~\cite[pp 200--211]{ito_mckean:diffusion_procs_sample_paths} for a proof of these facts using solutions to related differential equations. In this note we prove the following result on the almost sure rate of growth of $|N_t|$.

\begin{theorem}\label{thm:growth}
Suppose that the initial configuration consists of a finite number of particles at arbitrary positions in $\R$. Then, almost surely,
$$
\lim_{t\to\infty}\frac{\ln\ln|N_t|}{t}=2\sqrt{2\beta}.
$$
\end{theorem}

We define $R_t:=\max_{u\in N_t}Y_u(t)$ to be the right-most particle in the $(\beta y^2;\R)$-BBM. 
In Harris and Harris~\cite{harris_harris:inhombbm} it was shown that
\begin{equation}\label{eqn:rightmost_limit}
\lim_{t\to\infty}\frac{\ln R_t}{t}=\sqrt{2\beta}
\end{equation}
almost surely, and this result is crucial in our proof of Theorem~\ref{thm:growth} because it allows us both some control over the maximum breeding rate in the BBM, and also to show that the growth rate of $|N_t|$ is dominated by particles with spatial positions near $R_t$. Indeed, a BBM in which every particle has branching rate $\beta R_t^2$ would see its population grow as in Theorem~\ref{thm:growth},  and this point will provide our upper bound. Furthermore we shall see that, for any $\delta>0$, at all sufficiently large times $t$, a single particle located near $R_t$ will single-handedly build during time interval $[t,t+\delta]$ a progeny so large that it is again of the magnitude given by Theorem \ref{thm:growth}, which will prove the lower bound.

\section{Proof of the growth rate}

To prove Theorem~\ref{thm:growth}, it is sufficient to consider the case of the initial configuration of particles being a single particle at position $x\in \R$.

\medskip

\noindent\textbf{Upper bound.} Since the breeding potential is symmetric about the origin, we have from equation~\eqref{eqn:rightmost_limit} that for $\eps >0$ fixed, there exists, almost surely, a random time $\tau<\infty$ such that for all $u\in N_t$,
$$
\ln |Y_u(t)| < (\sqrt{2\beta}+\eps)t, \quad \text{for all }t>\tau.
$$
As a consequence of this, the breeding rate of any particle in the population is bounded by $\bar{\beta}_t:=\beta e^{2(\sqrt{2\beta}+\eps)t}$ for all $t>\tau$. 

We now introduce a coupled branching process $(\bar N_t, t\ge 0)$ as follows. For each $t \ge 0$ it consists in a population of particles $\{u : u \in \bar N_t\}$ whose
positions are denoted by $\{\bar Y_u (t) : u \in \bar N_t\}.$ Until time $\tau$ the two processes $N_t$ and $\bar N_t$ coincide and for all $t \ge \tau$ we will have $N_t \subseteq \bar N_t.$ Furthermore, we want that after time $\tau$ all particles in $\bar N_t$ branch at rate $\bar \beta_t.$

In order for this to make sense we must construct $(\bar N_t, t\ge 0)$ conditionally on  $( N_t, t\ge 0).$ Given $( N_t, t\ge 0),$ we let $( \bar N_t, t\le \tau)=( N_t, t\le \tau)$. For $t \ge \tau$, each particle $u \in N_t$ gives birth to an extra particle $v \in \bar N_t \backslash N_t$ at rate $\bar \beta_t - \beta Y_u(t)^2$ (note that if $u \in N_t$ then $Y_u(t) = \bar Y_u(t)$). Each particle thus created starts an independent BBM in $(\bar N_t, t\geq 0)$ with time-dependent branching rate $\bar \beta_t$ at time $t$. Thus it is clear that all particles in $\bar N_t$ branch at rate $\bar \beta_t$ for $t \ge \tau.$ 

It is furthermore clear that $$\vert N_t \vert \le  \vert \bar N_t \vert$$ for all $t\ge 0.$
The upshot is that after time $\tau$, $\bar N_t$ is a pure birth process, and as such is very well studied. 

Consider  $(Z_t, t\ge \tau)$  a pure birth process starting with a single particle at time $\tau$ with inhomogeneous rate $(\bar \beta_t, t\ge \tau)$, and define $$V_t := \inf \Big\{ s\geq\tau: \int_\tau^s \bar \beta_u\, \diffd u  \ge t\Big\}.$$ It is easily seen that $(Z_{V_t},t\geq 0)$ is simply a Yule process (i.e. a pure birth process where particles split into two at rate 1) and so we know that 
\[
 Z_{V_t} e^{-t} \to W_{\infty} \; \text{ as } t \to  \infty
 \]
 almost surely, where $W_{\infty}$ has an exponential mean 1 distribution. 
 However, since ${\int_\tau^{V_t} \bar \beta_u \,\diffd u}=t$ we see that this implies 
 \[
 Z_{t}\exp\Big( - \int_\tau^{t} \bar{\beta}_s\,\diffd s   \Big) \to W_{\infty} \; \text{ as } t \to  \infty.
 \]
 As $\bar N_t$ for $t\ge\tau$ is composed of the union of the
 offspring of the $|N_\tau|$ particles present at time $\tau$, we can write
$$
M_t := |\bar{N}_{t}|  \exp\Big( - \int_\tau^{t} \bar{\beta}_s\,\diffd s   \Big) 
=  \sum_{i=1}^{|N_{\tau}|} Z^{(i)}_t \exp\Big( - \int_\tau^{t} \bar{\beta}_s\,\diffd s   \Big)
$$
where the $Z^{(i)}$ are independent, identically distributed, copies of a pure birth process with birth rate
$\bar{\beta}_{t}.$ By conditioning on the value of $|N_\tau|$ (which is
almost surely finite), we see that
$$
M_t \to \bar{M}_\infty \in(0,\infty) \; \text{ as } t \to  \infty
$$
where, more precisely, $\bar{M}_\infty$ is distributed as the sum of $|N_\tau|$ independent exponential variables with mean 1.

Thus
$$
\ \left( \ln |\bar N_{t}| - \int_\tau^{t} \bar{\beta}_s\,\diffd s  \right) \to \ln \bar M_{\infty} \; \text{ as } t \to  \infty
$$ 
and since $\int_{\tau}^{t} \bar{\beta}_s\,\diffd s = \frac1c [ \exp(ct)-\exp(c\tau) ] $ with $c=2(\sqrt{2\beta}+\eps) $, it follows that
$$
\lim_{t\to\infty}\frac{\ln\ln |\bar{N}_{t}|}{t} = 2(\sqrt{2\beta}+ \eps),
$$
almost surely. Finally, since $\eps$ may be arbitrarily small,  using the coupling constructed above, we have that 
$$
\limsup_{t\to\infty}\frac{\ln\ln |N_t|}{t}\leq 2\sqrt{2\beta},
$$
almost surely, as required.

\medskip
\noindent
\textbf{Lower bound.} Fix $\eps>0$. From~\eqref{eqn:rightmost_limit}, we know that there exists, almost surely, a random time $\tau'<\infty$ such that for all $t>\tau'$, 
\begin{equation}\label{eqn:rightmost_bound}
R_t>e^{(\sqrt{2\beta}-\eps)t}.
\end{equation}

Fix $\delta >0$ and let $t_n:=  n\delta, n\in\{0,1,2,\ldots\}.$
We will show that it is sufficient to consider only the offspring of the rightmost particles $R_{t_n}$ to prove the lower bound for the global population growth rate. There are two steps to the argument. First, we will couple the sub-population descended from $R_{t_n}$ during the time interval $[t_n, t_{n+1}]$ to a BBM in which the expected population size remains finite. Second, we will show that, as $n$ tends to infinity, the population sizes of the coupled processes grow sufficiently quickly. For this, we use equation~\eqref{eqn:rightmost_bound} to give us a lower bound on the breeding rate in the coupled processes.

Let $N^{(n)}_s$ be the set of descendants of $R_{t_n}$ at time $s\in[t_n,t_{n+1}]$. As for the upper bound we introduce a coupled process $(\widetilde N^{(n)}_s, s\in [t_n, t_{n+1}]).$ This time we have that for all $s\in [t_n, t_{n+1}], \widetilde N^{(n)}_s \subseteq N^{(n)}_s \subseteq N_s.$ More  precisely, $(\widetilde N^{(n)}_s, s\in [t_n, t_{n+1}])$ is obtained from $( N^{(n)}_s, s\in [t_n, t_{n+1}])$ by cancelling some of the split events in $N^{(n)}_s,$ in the following way: if at time $s$ a particle $u \in  \widetilde N^{(n)}_s \subseteq N^{(n)}_s $ splits in the original process $N^{(n)}_s,$ it also splits in $\widetilde N^{(n)}_s$ 
with probability
$$
\frac{\widetilde{\beta}_n(Y_u(s))}{\beta Y_u(s)^2}\in[0,1],
$$
where
$$
\widetilde{\beta}_n(x):=
\begin{cases}
\beta e^{2(\sqrt{2\beta}-2\eps)t_n} \quad &  \text{if }  x \geq e^{(\sqrt{2\beta}-2\eps)t_n},\\
0&\text{otherwise.}
\end{cases}
$$
If the split event is rejected, one of the two offspring in $N^{(n)}_s$ is chosen at random to be the one which we keep in $\widetilde N^{(n)}_s.$

The process $(\widetilde N^{(n)}_s, t_n \le s \le t_{n+1})$ is thus a BBM started at time $t_n$ from a single particle at position $R_{t_n}$ with  space-dependent branching rate $\widetilde{\beta}_n$. Observe that we trivially have $\widetilde N^{(n)}_s \subseteq N^{(n)}_s$ for all $s\in[t_n,t_{n+1}],$ as announced, since we obtain $\widetilde N^{(n)}_s$ from $N^{(n)}_s$ merely by erasing some particles from $N^{(n)}.$

We also require another process on the same probability space, coupled to $(\widetilde{N}_s^{(n)}, s\in[t_n,t_{n+1}])$. This process is denoted $(\widehat{N}^{(n)}_s,s\in[t_n,t_{n+1}])$, and is defined by adding particles to $(\widetilde N^{(n)}_s,s\in[t_n,t_{n+1}])$ in such a way that every particle in $\widehat{N}^{(n)}$ breeds at constant rate $\beta e^{2(\sqrt{2\beta}-2\eps)t_n}$, irrespective of its spatial position.

More specifically, for $s\in[t_n,t_{n+1}]$, every particle $u\in \widetilde{N}_s^{(n)}$ gives birth to an extra particle $v\in\widehat{N}_s^{(n)}\backslash\widetilde{N}_s^{(n)}$ at rate 
$$
\beta e^{2(\sqrt{2\beta}-2\eps)t_n}\mathbf{1}_{\{Y_u(s)<e^{(\sqrt{2\beta}-2\eps)t_n}\}},
$$
and each particle thus created initiates an independent BBM in $(\widehat{N}^{(n)}_s,s\in[t_n,t_{n+1}])$ with constant breeding rate $\beta e^{2(\sqrt{2\beta}-2\eps)t_n}$. Note that we have $\widetilde N^{(n)}_s\subseteq\widehat{N}_s^{(n)}$ for all $s\in[t_n,t_{n+1}]$.

Now for $v \in N^{(n)}_{t_{n+1}}$ define the event
$$
A_n(v):=\Big\{\min_{s\in[t_n,t_{n+1}]}Y_v(s)<e^{(\sqrt{2\beta}-2\eps)t_n}   \Big\},
$$
and set $A_n=\cup_{v\in \widetilde{N}^{(n)}_{t_{n+1}}}A_n(v)$, i.e., the event that there exists a descendant $v$ of $R_{t_n}$ (in the modified process $\widetilde N^{(n)}$) such that $Y_v(s)<e^{(\sqrt{2\beta}-2\eps)t_n}$ for some $s\in [t_n, t_{n+1}]$. Observe that 
\begin{equation}\label{eqn:key_event}
\Big\{\widetilde{N}_{t_{n+1}}^{(n)}\neq\widehat{N}_{t_{n+1}}^{(n)}\Big\}\subseteq A_n.
\end{equation}

We also define the event $B_n$ as
$$
B_n:=\Big\{R_{t_n}>e^{(\sqrt{2\beta}-\eps)t_n}\Big\},
$$
and let $\{\F_t\}_{t\geq 0}$ be the natural filtration for the $(\beta y^2;\R)$-BBM. Our aim is to show that, almost surely, only finitely many of the events $A_n$ occur, from which it follows that the populations of the coupled subprocesses $\widetilde{N}^{(n)}$ are equal to the populations of the $\widehat{N}^{(n)}$ processes for all sufficiently large $n$, almost surely. Then the final step is to show that the populations $\widehat{N}^{(n)}$ grow sufficiently quickly to imply the desired lower bound on the size of the original population.

We start by writing the event of interest as
$$
A_n=(A_n\cap B_n)\cup (A_n\cap B_n^c),
$$
and we recall from equation~\eqref{eqn:rightmost_limit} that $P(\limsup B_n^c)=0$, and hence only finitely many of the events $A_n\cap B_n^c$ occur. We now use a standard `many-to-one' argument (see, for example, Hardy and Harris~\cite{hardy_harris:foundations}) to bound the probabilities $P(A_n\cap B_n|\F_{t_n})$.

Let $\Prob^x$ be the law of a driftless Brownian motion $Y$ started at the point $x\in\R$ (and $\E^x$ the expectation with respect to this law). Observing that $B_n\in\F_{t_n}$, we have
\begin{align*}
\mathbf{1}_{B_n}P(A_n | \F_{t_n})&\leq \mathbf{1}_{B_n}E \bigg[ \sum_{v\in \widetilde{N}^{(n)}_{t_{n+1}}}\mathbf{1}_{A_n(v)} \bigg| \F_{t_n}\bigg] \\
&= \mathbf{1}_{B_n}\expct^{R_{t_n}}\bigg( \exp\Big( \int_0^\delta \widetilde{\beta}_n(Y_s)\,\diffd s \Big);   \min_{ s\in [0,\delta]} Y_s < e^{(\sqrt{2\beta} -2\eps)t_n}\bigg)\\
&\leq\exp\Big({\beta e^{2(\sqrt{2\beta}-2\eps)t_n} \delta\Big)}\,\Prob^{e^{(\sqrt{2\beta}-\eps)t_n}}\left( \min_{ s\in [0,\delta]} Y_s < e^{(\sqrt{2\beta} -2\eps)t_n} \right).
\end{align*}
Using the reflection principle, we obtain that there exists $ C>0 $ such that 
\begin{align*}
\Prob^{e^{(\sqrt{2\beta}-\eps)t_n}} \left( \min_{ s\in [0,\delta]} Y_s < e^{(\sqrt{2\beta} -2\eps)t_n}  \right) 
&  = 2 \Prob^0  \Big(  Y_\delta < e^{(\sqrt{2\beta}-\eps)t_n}(e^{ -\eps t_n}  -1) \Big)\\
&  \le C \exp\Big(- \frac{1}{2\delta}(e^{(\sqrt{2\beta}-\eps)t_n}(e^{ -\eps t_n}  -1))^2\Big).
\end{align*}

Combining this series of inequalities gives, almost surely, a faster than exponentially decaying upper bound on $P(A_n\cap B_n|\F_{t_n})$, and so we have shown that
$$
\sum_{n\geq 0}P(A_n\cap B_n|\F_{t_n})<\infty
$$
almost surely. Since $A_n\cap B_n\in\F_{t_{n+1}}$, L\'{e}vy's extension of the Borel-Cantelli lemmas (see Williams~\cite[Theorem 12.15]{williams:probability_with_martingales}) lets us conclude that, almost surely, only finitely many of the events $A_n\cap B_n$ occur. Thus only finitely many of the events $A_n$ occur, almost surely.

Recalling equation~\eqref{eqn:key_event} we see that, almost surely, there exists a random integer $n_1 < \infty$ such that for all $n> n_1$, we have $\widetilde N^{(n)}_{s}=\widehat{N}_s^{(n)}$ for all $s\in[t_n,t_{n+1}]$.

The population size $|\widehat{N}_s^{(n)}|$ is a Yule process with constant breeding rate $\beta e^{2(\sqrt{2\beta}-2\eps)t_n}$, for $s\in[t_n,t_{n+1}]$. Hence $|\widehat{N}^{(n)}_{t_{n+1}}|$ has a geometric distribution with parameter 
$$
p_n:=\exp\Big( -\delta \beta e^{2(\sqrt{2\beta}-2\eps)t_n}\Big).
$$
If we define
$$
q_n:=\Big\lceil\exp\Big( e^{2(\sqrt{2\beta}-3\eps)t_n}\Big)\Big\rceil,
$$
then the probability that the number of particles in $\widehat{N}^{(n)}_{t_{n+1}}$ is smaller than $q_n$ is
$$
P(|\widehat{N}^{(n)}_{t_{n+1}}| \le q_n)=1-(1-p_n)^{q_n}=p_n q_n+o(p_nq_n).
$$
Using the Borel-Cantelli lemmas again, we see that there exists almost surely a random integer $n_2<\infty$ such that $|\widehat{N}^{(n)}_{t_{n+1}}| > q_n$ for all $n> n_2$.

Finally, for all $n>\max\{n_1, n_2\}$ and all $t\in[t_{n+1},t_{n+2}]$, we have that there are at least as many particles in $N_t$ as there are descendants of $R_{t_n}$ at time $t_{n+1}$, which is to say that $|N_t|\geq|\widetilde N^{(n)}_{t_{n+1}}|$. Hence, for $n$ sufficiently large,
\begin{equation}\label{eqn:final_ineq}
\frac{\ln\ln |N_t|}{t}\geq \frac{\ln\ln|\widetilde N^{(n)}_{t_{n+1}}|}{t_n}\frac{t_n}{t}\geq 2\sqrt{2\beta}-7	\eps,
\end{equation}
almost surely. (Certainly we must have $n>\max\{n_1, n_2\}$, but we may require that $n$ be larger still in order that the multiplicative factor $t_n/t$ is close enough to 1 for equation~\eqref{eqn:final_ineq} to hold.)

We can take $\eps$ to be arbitrarily small, and so obtain
$$
\liminf_{t\to\infty}\frac{\ln\ln |N_t|}{t}\geq 2\sqrt{2\beta},
$$
almost surely, which completes the proof.

\bibliographystyle{acm}
%\bibliography{../Resources/my_bibliography}

\end{document}